\documentclass[12pt]{article}

\usepackage{amsmath, amssymb, amsthm}
\usepackage{graphicx}
\usepackage{geometry}
\usepackage{mathtools}      
\usepackage{mathrsfs}       
\usepackage{bm}             
\usepackage{enumitem}       
\usepackage{thmtools}       
\usepackage{xcolor}         
\usepackage{hyperref}       
\usepackage{cleveref}       
\usepackage{tikz}           
\usepackage{caption}        
\usepackage{comment}
\usepackage[numbers]{natbib}

\newenvironment{keywords}{
  \noindent\textbf{Keywords:}\itshape
}{\par\vspace{1ex}}

\newenvironment{MSC}{
  \noindent\textbf{Mathematics Subject Classification (2020):}\itshape
}{\par\vspace{1ex}}

\newtheorem{theorem}{Theorem}[section]
\newtheorem{prop}[theorem]{Proposition}

\newtheorem{corollary}[theorem]{Corollary}

\newtheorem{conjecture}[theorem]{Conjecture}

\theoremstyle{definition}
\newtheorem{definition}[theorem]{Definition}

\newtheoremstyle{boldremark}
  {\topsep}   
  {\topsep}   
  {\normalfont}  
  {}          
  {\bfseries} 
  {.}         
  {.5em}      
  {}          

\theoremstyle{boldremark}
\newtheorem{remark}[theorem]{Remark}

\newenvironment{Proof}[1]{%
  \par\vspace{1ex}\noindent{\bf Proof{#1}:}\quad%
}{%
  \vspace{1.5ex}%
}
\def\Thanks#1{\gdef\thefootnote{\arabic{footnote}}\thanks{#1}}
\def\ThanksComma#1{\gdef\thefootnote{\arabic{footnote},}\thanks{#1}{
}}
\def\disp{\displaystyle}
\def\eps{\varepsilon}

\title{On Conjectures concerning the Labeled Coupon Collector Problem}
\author{D. Barak-Pelleg\ThanksComma{Department of Computer Science, Sami Shamoon College of Engineering, Beer Sheva 8410802, Israel.
E-mail: dina.barak.pelleg@gmail.com}
\Thanks{Research supported in part by a Hillel Gauchman scholarship.}
\and
D.~Berend\ThanksComma{Institute for the Theory of Computing and Department Mathematics, Ben-Gurion
University, Beer Sheva 84105, Israel.
E-mail: berend@bgu.ac.il}
\Thanks{Research supported in part by the Milken
Families Foundation Chair in Mathematics.}}
\date{}

\begin{document}

\maketitle

\begin{abstract}
We study a labeled variant of the classical Coupon Collector Problem (CCP), recently introduced by Tan et al., where coupons arrive in groups and only the set of labels is revealed. The goal is to determine the expected number of group drawings required to uniquely identify the labeling of all coupons. We focus on the case where groups consist of pairs ($k=2$), and provide rigorous proofs for two conjectures posed by Tan et al.
\end{abstract}

\begin{keywords}
Coupon Collector Problem, Labeling, Random Graphs
\end{keywords}

\begin{MSC}
Primary: 60C05; Secondary: 05C80, 68W20
\end{MSC}
\section{Introduction}

The classical \emph{Coupon Collector Problem} (CCP) is a well-known probabilistic model that has been extensively studied in the literature. In its standard form, one considers a collection of $n$ distinct coupons, each obtained independently and uniformly at random in successive trials. The central question is: how many trials are needed, on average, to collect all $n$ coupons? The expected number of trials required is known to be $n H_n\approx n(\log n+\gamma)$, where $H_n = \sum_{i=1}^n \frac{1}{i}$ is the $n$-th harmonic number and $\gamma$ is the Euler-Mascheroni constant. This result is foundational and appears in numerous probabilistic and algorithmic contexts (cf. \cite{erdos1961classical,flajolet1992birthday,motwani1995randomized}).

Variants of the CCP have been proposed to model more complex scenarios, including group arrivals, partial information, and labeling constraints. In this work, we consider the following version, studied by Tan et al \cite{tan2025labeledcoupon}. There coupons arrive in groups of $k$ distinct coupons each time, where $k$ is an arbitrary fixed integer between 1 and $n-1$. In each stage, each of the $\binom{n}{k}$ groups of size $k$ of coupons has the same probability of being drawn, independently of other stages. Each coupon $c$ has its own unique label $L(c)$, where $L$ is a bijection from the set of coupons to a set of size $n$ of so-called labels. When we get a group $\{c_1,c_2,\ldots,c_k\}$ of coupons, we are informed of the set $\{L(c_1),L(c_2),\ldots,L(c_k)\}$ of corresponding labels, but not which of the $k$ labels corresponds to each of the $k$ coupons.

The question is how many groups need to be drawn on average for us to be able to deduce uniquely what the label of each coupon is. In fact, the question has two variants. One is when the set of labels is known in advance, and the other when it is not. The difference is that, when the set of labels is known, it suffices to find out the labels of $n-1$ of the coupons; the label of the remaining coupon can then be inferred. The (random) number of required groups in the first case is denoted by $Q_I(n,k)$ and in the second by $Q_{II}(n,k)$. Tan et al \cite{tan2025labeledcoupon} considered the minimum number of groups enabling one to learn the function $L$ under both versions, as well as the expected values of $Q_I(n,k)$ and $Q_{II}(n,k)$. They raised two conjectures regarding the case $k=2$:
\begin{conjecture} \label{conj-1}
$E(Q_{II}(n,2))-E(Q_I(n,2))= \frac{1}{2} n + o(n)$.
\end{conjecture}

\begin{conjecture} \label{conj-2}
$E(Q_{II}(n,2))= \frac{1}{2} n H_n+ o(n)$, where $H_n=1+1/2+\cdots +1/n$ is the \( n \)-th harmonic number.
\end{conjecture}

In the sequel we prove both conjectures. The main results are formulated in Section \ref{sec:main_results} and the proofs are presented in Section \ref{sec:proofs}. We note that the error terms are much smaller than in the conjectures.

\section{Main Results} \label{sec:main_results}

\begin{theorem} \label{thm:EQI}
$E(Q_I(n,2))= \frac{1}{2} n H_n -\frac{1}{2} n + O(\log^5 n)$.
\end{theorem}

\begin{theorem} \label{thm:EQII}
$E(Q_{II}(n,2))= \frac{1}{2} n H_n+ O(\log^5 n)$.
\end{theorem}

\begin{corollary}
Both Conjectures \ref{conj-1} and \ref{conj-2} hold. Moreover, the error terms are $O(\log^5 n)$.
\end{corollary}

\begin{remark}
The error terms $O(\log^5 n)$ in Theorems \ref{thm:EQI} and \ref{thm:EQII} arise from our choice of the parameter $c = 2$ in the proofs (specifically in the application of Propositions {prop:conditional-completion} and {prop:conditional-almost-completion}). By choosing $c = 1 + \eps$ for arbitrarily small $\eps > 0$, these error terms can be improved to $O(\log^{3+o(1)}n)$. We have opted for $c = 2$ to keep the presentation simple and transparent.
\end{remark}

In the course of the proofs, we will employ the following result, which seems to be of independent interest. Note that, unlike the preceding theorems, it refers to the regular coupon collector -- one coupon at a time.

\begin{prop}
\label{prop:conditional-completion}
Let $T_k$ denote the additional time required to complete the coupon collection after obtaining exactly $k$ distinct coupons in the classical coupon collection problem, $0\le k\le n-2$. Let $c$ be a positive number, small enough so that $\log^c m < m$ for every integer $m\ge 2$.
Then, as $n\to\infty$:
\begin{equation} \label{bound-on-time-1}
    P\left(T_k \leq n\log(n-k) - cn\log\log(n-k)\right) \leq (1+o_n(1))\cdot \exp\left(-\log^c(n-k)\right).
\end{equation}
\end{prop}

\begin{remark}
\begin{enumerate}[label=\textbf{(\arabic*)}]
    \item \label{rem:part1} Similarly to the calculation of the expected time required to complete a full collection from the beginning, namely $E(T_0)$, we easily see that $E(T_k)=n H_{n-k}$. Unless $n-k$ is very small, this time is very close to $n(\log(n-k)+\gamma)$. It is interesting that relatively small deviations of $T_{n-k}$ downwards are much less likely than deviations upwards. For example, take $k=0$. It follows from the proposition that the probability for the collector to finish the collection $cn\log\log n$ steps ``too early'' for a fixed $c>1$ decays super-polynomially as a function of $n$. On the other hand, the probability of much larger upward deviation, namely $\alpha n\log n$ for any fixed $\alpha >0$, decays only polynomially \cite{chafai2010coupon}. 
    \item \label{rem:part2} If $k$ is very close to $n$, the right-hand side of \eqref{bound-on-time-1} does not decay to 0 as $n\to\infty$.
    \end{enumerate}
\end{remark}

\begin{prop}
\label{prop:conditional-almost-completion}
Let $T'_{k}$ denote the additional time required to get to $n-1$ collected coupons after obtaining exactly $k$ distinct coupons in the classical coupon collection problem, $0\le k\le n-2$. Let $c>0$ be small enough to ensure that $\log^c m < m$ for every integer $m\ge 2$. Then, as $n\to\infty$:
\begin{equation} \label{bound-on-time-2}
P\left(T'_{k} \leq n\log(n-k) - cn\log\log(n-k)\right) \leq 
(1+o_n(1))\cdot(n-k)\exp\left(-\log^c (n-k) \cdot \dfrac{n-k-1}{n-k}\right)
\end{equation}
\end{prop}

Since $E(T'_{k})=E(T_k)-n$, it makes sense to be interested in the probability of both $T_k$ and $T'_{k}$ being less than the same quantity $n\log(n-k) - cn\log\log(n-k)$. Unless $n-k$ is very small, a deviation of the order of magnitude $n\log\log (n-k)$ is much larger than the expectations difference. 

Theorems \ref{thm:EQI} and \ref{thm:EQII} relate to drawings of pairs. Thus, we will need not the variables $T_k$ and $T'_{k}$, but rather the variables denoting the analogous times when we draw pairs. As a first approximation, it is clear that the time required to complete a collection when drawing pairs is about half the time required when drawing single coupons. In fact, it tends to be even slightly less, because the second coupon in each pair is guaranteed to be different from the first. However, the advantage it creates is negligible, and both Propositions \ref{prop:conditional-completion} and \ref{prop:conditional-almost-completion} hold when we cut in half the time $n\log(n-k) - cn\log\log(n-k)$.

\section{Proofs}\label{sec:proofs}

\begin{Proof}{ of Theorem~\ref{thm:EQII}}

Consider the set $C$ of all coupons as the set of vertices of a graph $G=(C,E)$. In the beginning, $E$ is empty. At each step, if a pair of coupons $\{c_1,c_2\}$ is drawn, we add the edge $(c_1,c_2)$ to the graph (if it has not been added earlier). Clearly, at each stage we know the label of each vertex of degree 2 and above. Moreover, we know the label of each neighbor of a vertex with a known label. In other words, decomposing the graph into its connected components, we know the labels of all vertices belonging to components of size 3 and above, but do not know the labels of vertices in components of sizes 1 and 2.

Let $T$ denote the number of drawings of pairs of coupons until each coupon has been drawn at least once. In terms of the graph $G$ defined above, $T$ is the number of drawings until the graph has no isolated vertices, namely components of size 1. As explained above, $Q_{II}(n,2)$ is the number of drawings until all components are of size at least 3. In particular, $Q_{II}(n,2) \ge T$. By \cite[Thm. 2.9]{berend2025schilling},
\begin{equation} \label{time-with-pairs}
E(T)=\frac{1}{2} n H_n + O(\log n),
\end{equation}
so we only have to show that the expected number of additional drawings required to know all labels, after all coupons have been seen at least once, is $O(\log^5 n)$.

It will be convenient to view now the process of drawing the coupons in an additional, slightly different, way. Consider the coupons as arriving only one at a time, but each coupon that arrives at an even step $2m$ is drawn uniformly among all coupons except for the one that has just arrived at step $2m-1$. Denote the coupons according to their first arrivals -- coupon 1 is the first to be seen, coupon 2 -- the second, and so forth until the last coupon, that we denote by $n$, arrives.

Which pairs of coupons can form connected components of size 2 at the time we have seen all coupons at least once? Clearly, only pairs of the form $\{k,k+1\}$ for some $1\le k\le n-1$ can form such components. Specifically, if these two coupons arrive as a pair in the same drawing, and neither of them arrives again, except perhaps together with the other, until we finish seeing all coupons, then they will form a component. Note that the probability that $k$ and $k+1$ arrive together is $\frac{n-k}{n-1}$ if $k$ arrives at an odd step and $0$ otherwise. In particular, letting $I(k,k+1)$ denote the event that $k$ and $k+1$ arrive initially together, we have $P(I(k,k+1)) \le \frac{n-k}{n-1}$.

Suppose that some coupons $k,k+1$ arrived as a pair. We want to bound the probability that neither of the two will arrive again, together with one of the other $n-2$ coupons, until we see all coupons. According to Proposition \ref{prop:conditional-completion}, there is a probability of at least $1-(1+o(1))\cdot \exp\left(-\log^c(n-k)\right)$ that at least $\frac{1}{2}(n\log(n-k) - cn\log\log(n-k))$ pairs of coupons will be drawn until all coupons are seen. Here, $c$ may be any positive constant; for simplicity, we choose $c=2$. If this number of pairs is drawn, then the probability that none of the $2n-4$ pairs of the form $\{i,j\}$, with $i=k,k+1$ and $j\ne k,k+1$ is drawn in the process is
\begin{equation} \label{prob-stay-pair}
\begin{array}{rcl}
\left(1-\frac{2n-4}{n(n-1)/2}\right)^{\frac{1}{2}(n\log(n-k) - 2n\log\log(n-k))}
& = & \left(1-\frac{4}{n}+O(1/n^2)\right)^{\frac{1}{2}(n\log(n-k) - 2n\log\log(n-k))}\\\\
& = & \left(\left(1-\frac{4}{n}+O(1/n^2)\right)^{n/4}\right)^{2(\log(n-k) - 2\log\log(n-k))}\\\\
& = & \left(e^{-1}(1+O(1/n))\right)^{2(\log(n-k) - 2\log\log(n-k))}\\\\
& = & O\!\left(\frac{\log^{4}(n-k)}{(n-k)^2}\right).\\\\
\end{array}
\end{equation}

Denote by $F(k,k+1)$ the event that $\{k,k+1\}$ forms a connected component at the time all coupons have arrived. It follows from \eqref{prob-stay-pair} that
\begin{equation} \label{prob-stay-pair-2}
\begin{array}{rcl}
P(F(k,k+1) \mid I(k,k+1))
& = & \disp O\left(\frac{\log^{4}(n-k)}{(n-k)^2}\right) + (1+o(1))\cdot \exp\left(-\log^2(n-k)\right)\\\\
& = & \disp O\left(\frac{\log^{4}(n-k)}{(n-k)^2}\right).
\end{array}
\end{equation}
Consequently:
\begin{equation} \label{prob-be-component}
\begin{array}{rcl}
P(F(k,k+1)) & = & P(I(k,k+1)) \cdot P(F(k,k+1) \mid I(k,k+1))\\\\
& \le & \disp\frac{n-k}{n-1} \cdot O\left(\frac{\log^{4}(n-k)}{(n-k)^2}\right) \\\\
& = & \disp O\left(\frac{\log^{4}(n-k)}{n-k}\right) \cdot \frac{1}{n}.
\end{array}
\end{equation}
Let $X$ denote the number of connected components of $G$ when all coupons have arrived. Clearly,
$$E(X)=\sum_{k=1}^{n-1} P(F(k,k+1)),$$
and therefore
\begin{equation} \label{prob-be-component-2}
\begin{array}{rcl}
E(X) & = & \sum_{k=1}^{n-1}  
\disp O\left(\frac{\log^{4}(n-k)}{n-k}\right) \cdot \frac{1}{n}\\\\
& = & \disp \frac{1}{n} \cdot O\left(\int_1^n \frac{\log^4 x}{x} dx \right)\\\\
& = & \disp O\left(\frac{\log^5 n}{n}\right).
\end{array}
\end{equation}
Suppose $X>0$, so that we continue collecting pairs of coupons until no components of size~2 remain. The number of such components is reduced (by either 1 or 2) when one of the $2X$ coupons in the size-2 components arrives along with any coupons but its partner in the same component. It is easy to see that the probability for such a pair to arrive is at least $X/n$, and therefore one (or even two) of the size-2 components will be merged with another component in expected time at most $n/X$. Hence the expected time for all size-2 components to be merged with other components is bounded above as follows:
\begin{equation} \label{merge-time-conditioned}
\begin{array}{rcl}
E(Q_{II}(n,2)-T \mid X) & \le & 
\disp \frac{n}{X} + \frac{n}{X-1} + \cdots + \frac{n}{1}\\\\
& = & n H_X \\\\
& \le & nX.
\end{array}
\end{equation}
Thus:
\begin{equation} \label{merge-time}
E(Q_{II}(n,2)-T) \le 
E(nX)  = O\left(\log^5 n\right).
\end{equation}
In view of the initial discussion, this completes the proof.
\qed
\end{Proof}

\begin{Proof}{ of Theorem~\ref{thm:EQI}}
We proceed similarly to the proof of Theorem \ref{thm:EQII}, using also the same notations. 

Let $T'$ denote the number of drawings of pairs of coupons until $n-1$ of the coupons have been drawn at least once. In terms of the graph $G$ defined above, $T'$ is the number of drawings until the graph has at most one isolated vertex.

First, we want to find $E(T')$. Note first that, with a probability of at most $\frac{1}{n-1}$, coupons $n-1$ and $n$ arrive together as a pair, so we have drawn all coupons. Otherwise, if we continue the process of drawings (which is not required if there are no size-2 components at this point), the number of drawings until coupon $n$ is drawn is $G(2/n)$-distributed. It follows that
\begin{equation}
E(T) = E(T') + O(1/n)\cdot 0 + \left(1-O(1/n)\right)\cdot n/2.
\end{equation}
By \eqref{time-with-pairs}, this implies
\begin{equation} \label{time-for-n-th}
E(T') = \frac{1}{2} n H_n -\frac{1}{2}n + O(\log n).
\end{equation}
It thus remains to show that the expected number of drawings required to know all labels, after $n-1$ coupons have been seen at least once, is $O(\log^5 n)$.

We want to bound the probability that coupons $k,k+1$, that arrived as a pair, will not arrive again before $n-1$ arrives, except possibly again as a pair. According to Proposition \ref{prop:conditional-almost-completion}, there is a probability of at least 
$$1 - (1+o(1))\cdot(n-k)\exp\left(-\log^c (n-k) \cdot \dfrac{n-k-1}{n-k}\right)$$
that at least $\frac{1}{2}(n\log(n-k) - cn\log\log(n-k))$ pairs of coupons will be drawn until all coupons are seen. We take again $c=2$. If this number of pairs is drawn then, by \eqref{prob-stay-pair}, the probability that none of the $2n-4$ pairs of the form $\{i,j\}$, with $i=k,k+1$ and $j\ne k,k+1$ is drawn in the process is
$O\left(\frac{\log^{4}(n-k)}{(n-k)^2}\right).$

Let $I(k,k+1)$ be as in the proof of Theorem \ref{thm:EQII}, and $A'_F(k,k+1)$ be the same as $F(k,k+1)$, except that it refers to the state after coupon $n-1$ has arrived. Similarly to~\eqref{prob-stay-pair-2}, we now have 
\begin{equation} \label{prob-stay-pair-3}
\begin{array}{rcl}
P(A'_F(k,k+1) \mid I(k,k+1))
& = & \disp O\left(\frac{\log^{4}(n-k)}{(n-k)^2}\right)\\\\
& & + (1+o(1))\cdot \exp\left(-\log^2(n-k)\cdot \dfrac{n-k-1}{n-k}\right)\\\\
& = & \disp O\left(\frac{\log^{4}(n-k)}{(n-k)^2}\right).
\end{array}
\end{equation}
The rest of the proof stays intact with the following minor changes:
\begin{itemize}
\item $X$ is replaced by $X'$ -- the number size-2 components when coupon $n-1$ arrives.
\item $F(k,k+1)$ is replaced by $A'_F(k,k+1)$.
\item $Q_{II}(n,2)$ is replaced by $Q_I(n,2)$.
\end{itemize}
\qed
\end{Proof}

The proof of Proposition \ref{prop:conditional-completion} will require an auxiliary result, for which we introduce the following notion.
\begin{definition}\cite{dubhashi2009concentration,joag1983negative}
\label{def:negative-association}
Random variables $X_1, X_2, \ldots, X_n$ are \emph{negatively associated} if for every disjoint index sets $I, J \subseteq [n]$ and coordinate-wise non-decreasing functions $f$ and $g$, we have
\begin{equation}
    \text{Cov}(f(X_I), g(X_J)) \leq 0,
\end{equation}
where $X_I = (X_i : i \in I)$ and $X_J = (X_j : j \in J)$.
\end{definition}

This notion is relevant to us through the balls and bins problem. Suppose $m$ balls are tossed independently into $n$ bins, each ball having probability $1/n$ of landing at any bin. Then the numbers of balls landing in the various bins are negatively associated \cite{dubhashi2009concentration,joag1983negative}. Viewing coupon types as bins, we see that the numbers of coupons of each type obtained in a series of drawings are negatively associated.

An important property of negatively associated variables $X_1, X_2, \ldots, X_n$ is that
$$E\left(\prod_{i=1}^n f_i(X_i)\right) \le \prod_{i=1}^n E(f_i(X_i))$$
for every non-decreasing functions $f_1,\ldots,f_n$. This implies that, if $I$ is a set of coupons, and $A_i$ denotes the event that coupon $i\in I$ has been obtained in some sequence of drawings, then
$$P\left(\bigcap_{i\in I} A_i\right) \le \prod_{i\in I}P(A_i).$$

\begin{Proof}{ of Proposition~\ref{prop:conditional-completion}}

After collecting $k$ distinct coupons, we need to collect the remaining $n-k$ coupons. Let $B_j, 1\le j\le n-k$, denote the event that the $j$-th remaining coupon is collected within $n\log(n-k)-cn\log\log(n-k)$ drawings. The probability that this coupon is \emph{not} collected in these trials is:
\begin{equation}
\begin{array}{rcl}
    P\left(\overline{B_j}\right) 
    & = & \left(1 - \frac{1}{n}\right)^{n\log(n-k)-cn\log\log(n-k)}\\\\
    & = & \exp\left(n\cdot (\log(n-k)-c\log\log(n-k)) \cdot \log(1-1/n) \right)\\\\
    & = & \exp\left(-\log(n-k)+c\log\log(n-k)+O\left(\frac{\log n}{n} \right)\right)\\\\
    & = & \dfrac{\log^c (n-k)}{n-k} \cdot\left(1+O\left(\frac{\log n}{n} \right)\right).
\end{array}
\end{equation}
Hence:
$$P(B_j) = 1 - \dfrac{\log^c (n-k)}{n-k} \cdot\left(1+O\left(\frac{\log n}{n} \right)\right).$$
As explained before the proof, the events $B_j$ for the $n-k$ remaining coupons are negatively associated. Therefore:
\begin{equation}
\begin{array}{rcl}
    \disp P\left(\bigcap_{j=1}^{n-k} B_j\right) 
    & \leq & \disp\prod_{j=1}^{n-k} P\left(B_j\right)\\\\
    & = & \disp\left(1 - \dfrac{\log^c (n-k)}{n-k} \cdot\left(1+O\left(\frac{\log n}{n} \right)\right)\right) ^{n-k}\\\\
    & \le & \exp\left(-\log^c (n-k) \cdot\left(1+O\left(\frac{\log n}{n} \right)\right)\right)\\\\
    & = & \exp\left(-\log^c (n-k)\right) \cdot\left(1+O\left(\frac{\log^{c+1} n}{n} \right)\right).\\\\
\end{array}
\end{equation}
Now the event considered in the proposition, namely $\{T_k \leq n\log(n-k) - cn\log\log(n-k)\}$, is exactly the intersection of all $B_j$'s, and therefore:
\begin{equation}
    P\left(T_k \leq n\log(n-k) - cn\log\log(n-k)\right)
    \leq (1+o(1))\cdot \exp\left(-\log^c(n-k)\right).
\end{equation}

\end{Proof}

\begin{Proof}{ of Proposition~\ref{prop:conditional-almost-completion}}
We proceed similarly to the proof of Proposition \ref{prop:conditional-completion}.    
Clearly,
\begin{equation} \label{bound-prob}
\{T'_{k}\le n\log(n-k) - cn\log\log(n-k)\}=\bigcup_{i=1}^{n-k}\left(\bigcap_{j\ne i}B_j\right),
\end{equation}
and hence, by the union bound and symmetry,
\begin{equation} \label{bound-prob-2}
P(T'_{k}\le n\log(n-k) - cn\log\log(n-k))
\le
(n-k)\,P\!\left(\bigcap_{j=1}^{n-k-1}B_j\right).
\end{equation}
Similarly to the proof of Proposition~\ref{prop:conditional-completion},
\begin{equation} \label{bound-intersection-prob}
\begin{array}{rcl}
\disp P\!\left(\bigcap_{j=1}^{n-k-1}B_j\right)
& \le &
\disp\prod_{j=1}^{n-k-1}P(B_j)\\\\
& = & \disp\left(1 - \dfrac{\log^c (n-k)}{n-k} \cdot\left(1+O\left(\frac{\log n}{n} \right)\right)\right) ^{n-k-1}\\\\
& = & \disp\left(1 - \dfrac{\log^c (n-k)}{n-k-1} \cdot \dfrac{n-k-1}{n-k} \cdot\left(1+O\left(\frac{\log n}{n} \right)\right)\right) ^{n-k-1}\\\\
& \le & \exp\left(-\log^c (n-k) \cdot \dfrac{n-k-1}{n-k} \cdot\left(1+O\left(\frac{\log n}{n} \right)\right)\right)\\\\
    & = & \exp\left(-\log^c (n-k) \cdot \dfrac{n-k-1}{n-k}\right) \cdot\left(1+O\left(\frac{\log^{c+1} n}{n} \right)\right).\\\\
\end{array}
\end{equation}
Finally, by \eqref{bound-prob-2} and \eqref{bound-intersection-prob},
\begin{equation} \label{final-bound}
\begin{array}{rcl}
P(T'_{k}\le n\log(n-k) - cn\log\log(n-k))
& \le & (n-k)\exp\left(-\log^c (n-k) \cdot \dfrac{n-k-1}{n-k}\right)\\\\
& &\cdot\left(1+O\left(\frac{\log^{c+1} n}{n} \right)\right).\\\\
\end{array}
\end{equation}
 \qed
\end{Proof}

\section*{Acknowledgements}
We would like to thank J.-P. Allouche for drawing our attention to \cite{tan2025labeledcoupon}, the starting point for this work.

\bibliographystyle{plainnat}
\bibliography{refs}

\end{document}